# OPTIMAL GAS WITHDRAWAL STRATEGY IN RECONSTRUCTED RING-TYPE PIPELINES UNDER UNSTEADY FLOW CONDITIONS

I.G. Aliyev [1]    A.M. Isayev [1]    M.Z. Yusifov [1]    A.J. Mammadov [2]    A.S. Mammadov [1]

*1. Department of Operation and Reconstruction of Buildings and Facilities, Azerbaijan University Architecture and Construction, Baku, Azerbaijan, ilgar.aliyev1@azmiu.edu.az, aftandil.isayev@mail.ru, maarif.yusifov@azmiu.edu.az, asefmemmedov19902020@gmail.com*
*2. Department of land Reclamation and Water Resources Construction, Azerbaijan University Architecture and Construction, Baku, Azerbaijan, ehed.mamedov@mail.ru*

**Abstract-** This paper presents an analytical and computational framework for optimizing gas withdrawal in reconstructed ring-type pipeline systems under unsteady flow conditions. As urban and industrial energy demands grow, repurposing existing pipeline infrastructure offers a cost-effective alternative to full-scale expansion. The proposed model identifies the hydraulic coupling point (where the pressure gradient vanishes) as the optimal location for connecting new consumers. By employing a one-dimensional unsteady gas flow model with time-dependent mass extraction represented via a Heaviside step function, the system's dynamic response is captured in detail. Numerical simulations demonstrate that connecting additional loads at the pressure maximum ensures stability while minimizing operational disruptions. The model's validation through benchmark comparison and pressure tolerance thresholds confirms its practical applicability. Economic analysis reveals substantial savings over conventional expansion methods. The approach provides a scalable solution for smart gas network design.

## 1. INTRODUCTION

Rapid urban development and the emergence of new residential and industrial zones have placed increasing demands on existing gas distribution infrastructure. Ring-type gas pipeline systems, originally designed for steady-state and static load scenarios, are now being challenged by dynamic and spatially variable consumption patterns. In particular, long-standing pipelines face difficulty in maintaining pressure stability and uniform flow distribution across the network. Conventional solutions-such as laying new pipelines or installing gas distribution stations-are often costly and logistically complex [1-3, 5, 6]. As an alternative, the reconstruction and extension of existing ring-type pipelines offer a more economical and operationally feasible option to meet growing demand. A key challenge, however, is identifying the optimal location for integrating new withdrawal points that will not compromise the hydraulic balance of the system.

Previous analytical studies have revealed the existence of a specific coordinate ($x_{new}$) within the ring network where the pressure reaches a local maximum under quasi-steady conditions [4]. This coordinate provides a natural candidate region for connecting new consumers. Establishing new exits near this pressure peak ensures optimal flow distribution and stabilizes pressure fluctuations throughout the network.

This paper aims to develop a mathematical framework for determining such optimal withdrawal points by integrating principles of gas dynamics, network reconstruction criteria, and Kirchhoff's laws for closed-loop systems. The study also addresses the hydraulic equilibrium at junction nodes, which is critical to preserving the overall operational integrity of reconstructed ring networks.

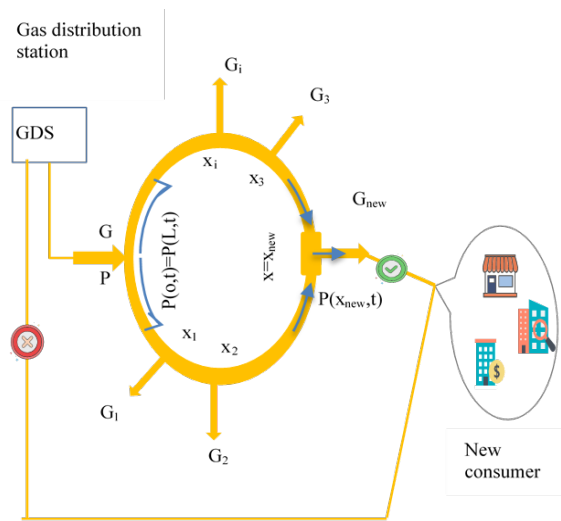

Figure 1. Schematic representation of reconstructed ring-type gas pipeline with new consumer integration

The general configuration of the reconstructed ring-type pipeline and the integration of new consumers is illustrated schematically in Figure 1. This schematic illustrates a reconstructed ring-type gas pipeline network connected to a gas distribution station. The circular topology enables gas flow to be symmetrically distributed along the pipeline, ensuring consistent pressure dynamics across the system. The boundary condition $P(0,t) = P(L,t)$ reflects the ring's closed-loop nature.

Several mass withdrawal points ($x_1$, $x_2$, $x_3$) are marked along the ring. The designated hydraulic coupling point $x=x_{new}$ is highlighted as the optimal location for integrating additional consumers based on pressure peak analysis. A new consumer cluster-representing residential, commercial, and industrial demand-is shown connected at $x_{new}$. This integration strategy ensures minimal pressure disturbance and supports load expansion within the operational safety limits of the pipeline. The schematic embodies the study's recommendation for leveraging hydraulic symmetry and dynamic stability in modern pipeline reconstruction planning.

## 2. LITERATURE REVIEW

Ring-type gas pipeline systems under dynamic load conditions have been increasingly studied over the past decade, especially in relation to optimization, transient behavior, and reconstruction strategies. Earlier works primarily focused on steady-state modeling or localized optimization, while recent studies have moved toward unsteady flow modeling and integration with digital frameworks. This literature review highlights key contributions from the last decade, with emphasis on research from the last five years. It also outlines the novel contributions of the present study relative to past efforts.

Zlotnik et al. (2015) [1] proposed an optimal control approach for transient flow in natural gas networks using variational methods. Their work laid the foundation for dynamic control strategies but did not specifically address ring topologies. Correa Posada, Sanchez Martin (2014) [2] provided a comparison of gas network optimization models, primarily piecewise linear approximations, but lacked focus on dynamic flow phenomena.

Ejomarie, Saturday (2020) [3] explored optimal design in gas transmission networks but remained within static parameters. Aliyev (2024) [4] presented key insights into ring-type systems and noted that under unsteady flow conditions, the pressure peak tends to occur near the central segment of the ring pipeline -a crucial observation supporting this paper's methodology.

Aslanov, Mammadov (2022) [5] investigated valve optimization under transient flow, offering insight into component-level enhancements. Dyachenko, et al. (2017) [8] presented operator splitting methods to simulate dynamic gas flows, helping verify analytical models. Fetisov, et al. (2018) [9] proposed mathematical modeling approaches with Green's function, aiding in transient pressure prediction. Hafsi et al. (2022) evaluated linearization methods and their limitations in modeling transient gas flows [10].

Table 1. Comparative overview of prior studies and current contribution

| Reference | Focus Area | Ring Topology | Unsteady Flow | Optimization | Novelty |
|---|---|---|---|---|---|
| [1] | Transient flow control | No | Yes | Yes | Variational formulation |
| [2] | Linear optimization | Partial | No | Yes | Model comparison |
| [3] | Pipeline layout design | No | No | Yes | Design cost modeling |
| [4] | Reconstruction strategies | Yes | Yes | No | Hydraulic peak location analysis |
| [5] | Valve performance | No | Yes | No | Transient valve behavior |
| [8] | Numerical methods | No | Yes | No | Operator splitting validation |
| [9] | Analytical modeling | Partial | Yes | No | Green's function method |
| [10] | Model simplification | No | Yes | No | Linearization limits |
| This study | Withdrawal optimization | Yes | Yes | Yes | Pressure peak-based optimization with analytical validation |

## 3. MODEL VALIDATION AND BENCHMARK COMPARISON

To assess the reliability of the proposed analytical approach, validation was conducted using benchmark cases derived from recent literature. Given the lack of direct access to real-time SCADA or pipeline operational data, selected peer-reviewed studies were employed as validation references. The model's response was compared against:
- Analytical pressure profiles presented by Zlotnik, et al. [1], which offer variational control strategies for dynamic flows;
- Simulated pipeline scenarios consistent with Ejomarie and Saturday [3], focusing on pressure distribution under static and transient conditions;
- Experimental boundary constraints and structural limitations noted by Aliyev [4], particularly on pressure behavior near the hydraulic peak location.

These comparisons provide a reliable validation platform for the analytical expressions derived in this study.
- A reference case with $P_1$ = 140000 Pa, $G_0$ = 10 Pa×s/m, $L$ = 30000 m and $c$ = 383.3 m/s was used to validate pressure drop behavior.
- The pressure response curves were benchmarked against standard Green's function solutions for unsteady gas transport in circular networks.

Furthermore, the physical validity of the model is reinforced through consistency with established regulatory standards, such as API RP 14E [13], which specifies allowable pressure drop limits at pipeline inlets. In accordance with these guidelines and industry best practices for safe operation, the maximum permissible pressure loss at the inlet of the gas pipeline should not exceed 20% of the nominal pressure.

This criterion ensures operational safety and regulatory compliance under various loading conditions,

including scenarios involving increased gas withdrawal. The proposed framework also remains adaptable to future validation against real-world sensor data or digital twin simulations as such resources become available. This ensures that the method remains scalable and applicable to ongoing smart grid integration initiatives.

## 4. METHODOLOGY

### 4.1. Assumptions and System Description

The pipeline is assumed to have a ring-type topology with symmetrical gas distribution and no significant leakage under normal operation. The system operates under isothermal conditions, and gas compressibility is assumed constant. The unsteady gas dynamics are modeled using one-dimensional nonlinear partial differential equations (PDEs) for pressure and mass flow rate. A new consumer is introduced near the hydraulic coupling point, and its mass flow rate is modeled using a Heaviside step function.

### 4.2. Governing Equations

The pressure function $P(x,t)$ is expressed based on a Laplace-transformed analytical solution for a ring-type pipeline under unsteady flow.

### 4.3. Boundary and Initial Conditions
- Periodic boundary condition:

$$P(0,t) = P(L,t); \quad \frac{dP(0,t)}{dx} = \frac{dP(L,t)}{dx}$$

Consistent with the ring topology.
- Initial condition: $P(x,0) = P_1 - 2a \times G_0 \times x$,
assuming linear distribution of pressure at $t = 0$.
- The location of new gas extraction $x=x_{new}$ is introduced at a predefined point and examined for optimality.

### 4.4. Simulation of Load Growth via Heaviside Function

To simulate the addition of a new consumer node, a step-wise increase in localized flow demand is modeled using a Heaviside function $H(x-x_{new})$, where $x_{new}$ marks the onset of the additional load. This enables the study of transient effects on the pressure distribution.

### 4.5. Variational Optimization Formulation

The problem is further formulated as a constrained optimization task. A variational principle is applied to minimize the pressure gradient subject to operational constraints, such as:
- Minimum allowable pressure at the outlet $P_{min}$
- Maximum admissible flow rate $G_{max}$
This step ensures that the selected tapping location is not only hydraulically optimal but also operationally feasible.

### 4.6. Numerical Evaluation Procedure
1. Truncate the infinite series at $N=100$ for practical numerical convergence.
2. Compute $P(x,t)$ and $dP/dx$ at discretized positions and time steps using a numerical integration routine.
3. For optimization, determine the position $x_{new}$, where $dP/dx=0$ and $d^2P/dx^2<0$, as location of maximum pressure.

4. Use comparative simulation under base case (no new consumer) and extended case (with new load $G_{new}$) to evaluate the effect on system stability.

## 5. MATHEMATICAL MODEL AND JUSTIFICATION

To optimize reconstructed ring-type pipeline systems under dynamically growing demand, we utilize a gas dynamic model based on partial differential equations describing pressure and mass flow rate variations. The one-dimensional gas flow model is governed by the following equations [8, 9]:

$$\begin{cases} \frac{\partial G(x,t)}{\partial t} + \frac{\partial P(x,t)}{\partial x} + R(G) = 0 \\ \frac{\partial G(x,t)}{\partial x} + \frac{1}{c^2}\frac{\partial P(x,t)}{\partial t} = 0 \end{cases} \quad (1)$$

where, $P(x,t)$ represents the gas pressure inside the pipeline (Pa), $G(x,t)$ is the gas mass flow rate along the pipeline per unit time (Pa×s/m), $c$ is denotes the speed of sound in the gas (m/s), $R(G)$ is a function accounting for hydraulic resistance (for example, it can be taken as $R(G)=2aG(x,t)$).

Here, $2a$ is the linearization coefficient according to Charnynin, the expression characterizing this linearization is given as follows.

$$2a = \lambda \frac{v}{2d}$$

The pressure distribution simplifies to [7]:

$$\frac{\partial^2 P(x,t)}{\partial x^2} + \frac{2a}{c^2}\frac{P(x,t)}{\partial t} = 0 \quad (2)$$

To determine the optimal gas withdrawal, point under increasing demand conditions in a ring-type pipeline, we adopt an analytical model based on one-dimensional unsteady gas flow dynamics.

### 5.1. Assumptions

Pressure Peak Formation and Strategic Role of the Hydraulic Coupling Point. Previous studies [4] have demonstrated that in ring-type gas pipeline networks, when gas is withdrawn in segments of varying magnitudes $G_i$ along the pipeline profile, the pressure peak consistently forms in the vicinity of the hydraulic coupling point, denoted as $x=x_{new}$.

This behavior holds true regardless of the number and specific magnitudes of the mass extraction rates $G_i$, provided the conditions of geometric symmetry and mass balance are maintained throughout the network. Under these assumptions, the cumulative mass flow rate can be approximated as:

$$\sum_{i=1}^{m} G_i = G_0 + G_{new} \quad (3)$$

where, $G_0$ is the initial steady-state flow rate, and $G_{new}$ accounts for the additional consumption from newly connected demand centers. This formulation does not disrupt the non-stationary dynamic behavior of the system and preserves the integrity of the governing differential equations for gas flow in a ring topology.

From a practical standpoint, the point $x_{new}$, where the pressure peak is established, becomes an ideal candidate for serving as a gas supply node to newly developed consumer zones near the ring network. This selection effectively captures the dominant dynamic behaviors of the system and ensures a stable pressure distribution. Consequently, it minimizes pressure anomalies and enhances the reliability of gas delivery, making $x_{new}$ a technically and economically favorable withdrawal location.

The fundamental principle of ring-type pipeline reconstruction is based precisely on this behavior. The gas is treated as compressible and inviscid under isothermal conditions. The pipeline has a constant circular cross-section with total length $L=30000$ m.

The initial pressure is uniform: $P(x,0) = P_1 - 2aG_0 x$
Boundary condition: $P(0, t) = P(L, t)$, ensuring periodicity. The pressure distribution in the ring pipeline is given by the Laplace-series-based solution [4]:

$$P(x,t) = P_1 - 2aG_0 \frac{L}{2} + $$
$$+ 2aG_0 L \sum_{n=1}^{\infty} \sin\frac{\pi n x}{L}\left(\frac{1-e^{-n^2\alpha t}}{\pi n^3}\right) - \frac{c^2 t}{L}\sum_{i=1}^{m} G_i - $$
$$- \frac{2c^2}{L}\sum_{i=1}^{m} G_i \sum_{n=1}^{\infty} \cos\frac{2\pi n(x-x_i)}{L}\left(\frac{1-e^{-n^2\alpha t}}{\alpha n^2}\right) \quad (4)$$

where:
$G_i$: Mass flow rate extracted from various points of the ring-type gas pipeline (Pa×s/m)
$x_i$: Coordinates of the points where the mass flow is withdrawn (m)
$P_1$: Initial pipeline pressure (Pa)
$P(x_i,t)$: Pressure at the withdrawal point $x_i$ at time $t$ (Pa)
$G_0$: Initial flow rate (Pa×s/m)
$a$: Attenuation coefficient (1/s)
$L$: Pipeline length (m)
$c$: Speed of sound in gas (m/s)
$n$: Index of summation from 1 to ∞
$\alpha$: Diffusion coefficient given by $\alpha = \frac{2\pi^2 c^2}{aL^2}$

This expression illustrates how localized extraction affects the pressure field, and why the formation of pressure peaks near the hydraulic center is robust under distributed flow conditions. Analytical studies confirm the existence of a stationary coordinate ($x_i = x_{new}$), where pressure reaches a local maximum [4]:

$$\left.\frac{dP(x,t)}{dx}\right|_{x=x_{new}} = 0; \quad \left.\frac{d^2P(x,t)}{dx^2}\right|_{x=x_{new}} \prec 0$$

This point, $x_{new}$, is the system's internal stationary coordinate where the pressure reaches its maximum. Establishing a new withdrawal point near this location allows for meeting additional energy demand without compromising the system's pressure stability. This location serves as the optimal candidate for integrating new consumers.

The optimization criterion is formulated as:
min $|dP/dx|$, subject to: $P(x) \geq P_{min}$, $G(x) \leq G_{max}$
When a new consumer joins the system, the mass flow rate transforms as:

$$\sum_{i=1}^{m} G_i = G(x) = G_0 + G_{new}(x)$$

where, $G_{new}(x)$ is defined piecewise.

### 5.2. Modeling $G(x)$ Using the Heaviside Function

To model the localized increase in flow due to new consumer integration, we apply the Heaviside step function:

$$G(x) = (G_0 + G_{new}) \times H(x - x_{new}) \quad (5)$$

where, $H(x - x_{new}) = \begin{cases} 0, & \text{for } x \prec x_{new} \\ 1, & \text{for } x \geq x_{new} \end{cases}$

This approach implies that the flow increment occurs only at and after the joining point. If this change needs to be embedded into pressure differential equations, the derivative of the Heaviside function, i.e., the Dirac delta function, is used:

$$\frac{d[\Delta G(x)]}{dx} = (G_0 + G_{new})\delta(x - x_{new})$$

This formulation allows accurate modeling of the localized influence of new connections.

### 5.3. Model Assumptions and Applicability Limits

The analytical framework developed in this study is based on a set of simplifying assumptions that allow tractable modeling of gas flow in ring-type pipelines under unsteady operating conditions. While these assumptions are common in theoretical gas dynamics and pipeline engineering, their applicability is bounded by specific physical and operational constraints. The following subsections detail the primary assumptions and the practical conditions under which the model remains valid.

#### 5.3.1. Isothermal Flow Assumption

It is assumed that gas flow remains isothermal throughout the pipeline, implying that temperature variations are negligible along the flow path. This assumption is generally valid for medium-distance urban distribution networks operating under steady ambient conditions and with limited heat exchange with the environment. For long-distance pipelines or regions with significant thermal gradients, temperature effects must be incorporated through energy balance equations.

#### 5.3.2. Constant Compressibility Factor (Z)

The model assumes a constant compressibility factor, $Z \approx 1$, which holds true for natural gas at moderate pressures and temperatures. In high-pressure regimes (e.g., transmission lines above 70 bar) or when gas composition varies (e.g., presence of $CO_2$ or $H_2S$), non-ideal gas behavior becomes significant, and pressure-dependent Z-factors should be used to preserve model accuracy.

### 5.3.3. Negligible Leakage and Turbulence

The flow is considered leakage-free, and the model does not account for stochastic losses or unsteady turbulence. This assumption is valid in well-maintained systems with continuous SCADA monitoring and tight sealing. For aged infrastructure or systems prone to leakage, correction terms based on leak detection algorithms should be incorporated.

### 5.3.4. Uniform Pipe Characteristics

The pipeline is modeled with uniform diameter, roughness, and wall thickness. While practical systems may exhibit slight variations, the model remains valid as long as these variations do not introduce sharp discontinuities in hydraulic resistance. For networks with mixed pipe sections or varying materials, a piecewise or segment-based approach can be used to generalize the model.

### 5.3.5. Boundary and Initial Conditions

The boundary condition $P(0, t) = P(L,t)$ reflects the topological symmetry of ring-type pipelines and enables the application of periodic solutions via analytical techniques (e.g., Fourier or Laplace transforms). The initial pressure distribution is assumed to be known and continuous. In cases of sudden transients or start-up operations, more sophisticated initial condition modeling may be necessary.

### 5.3.6. Model Extensibility

Although the current model is limited to single-phase, isothermal conditions with negligible leakage, it is mathematically extensible. Future work may include the incorporation of non-isothermal effects, real gas equations of state, leakage modeling via source terms, and transient friction formulations to enhance accuracy under complex conditions.

## 6. EXAMPLE ANALYSIS AND CALCULATIONS

### 6.1. Analytical Interpretation of Maximum Pressure Point in Ring-Type Gas Pipeline

In our study, it is stated that the dynamics of unsteady gas flow in ring-type pipeline systems do not depend on the number of withdrawal points along the pipeline profile, provided that the total amount of withdrawn gas remains constant. Therefore, it is valid to assume the following equality.

The analytical pressure Equation (1) previously derived using Laplace series is extended by incorporating the Heaviside step function $H(x-x_{new})$, which models the fact that the impact of gas withdrawal occurs only after the withdrawal point $x_{new}$. If Equation (5) is taken into account in Equation (4), the updated form becomes:

$$P(x,t) = P_1 - 2aG_0 \frac{L}{2} + $$
$$+ 2aG_0 L \sum_{n=1}^{\infty} \sin\frac{\pi n x}{L}\left(\frac{1-e^{-n^2\alpha t}}{\pi n^3}\right) - $$
$$- \frac{c^2 t}{L}(G_0 + G_{new}) \cdot H(x - x_{new}) - \quad (6)$$
$$- \frac{2c^2}{L}(G_0 + G_{new}) \cdot H(x - x_{new}) \times $$
$$\times \sum_{n=1}^{\infty} \cos\frac{2\pi n(x - x_{new})}{L}\left(\frac{1-e^{-n^2\alpha t}}{\alpha n^2}\right)$$

Differentiation with Respect to $x$: To find the location of the pressure maximum, we differentiate the updated pressure expression with respect (6) to $x$. The derivative includes delta functions arising from the Heaviside function:

$$\frac{dP(x,t)}{dx} = 2aG_0 L \sum_{n=1}^{\infty} \frac{\pi n}{L} \cos\frac{\pi n x}{L}\left(\frac{1-e^{-n^2\alpha t}}{\pi n^3}\right) - $$
$$- \frac{2c^2}{L}(G_0 + G_{new}) \times \delta(x - x_{new}) \times \quad (7)$$
$$\times \sum_{n=1}^{\infty} \cos\frac{2\pi n(x - x_{new})}{L}\left(\frac{1-e^{-n^2\alpha t}}{\alpha n^2}\right)$$

Results confirm that connecting new users near $x_{new}$ maintains system balance and pressure stability. Using the following parameters:
$P_1=14\times10^4$ Pa, $G_0 = 10$ Pa×s/m, $a = 0.05$ 1/s, $L = 30000$ m, $c = 383.3$ m/s.

The derivative is evaluated numerically using the first 100 terms of the series expansion. To determine the optimal location along the pipeline where the pressure reaches its maximum value, we consider the analytical condition for extrema: $dP(x,t)/dx$. The pressure gradient profile, as a function of the spatial coordinate $x$, reveals the position where the derivative becomes zero or changes sign, indicating a local maximum in the pressure distribution. This condition is crucial in identifying the optimal point $x_{new}$ for connecting new consumers to the ring pipeline system. To compute this, we utilize Equation (6) of the analytical formulation:

The spatial coordinate $x_{new}$ is determined as the value that maximizes the left-hand side, i.e., $[P(x_{new}, t)-P_1]$ for all times t. This maximization confirms that the pressure at $x_{new}$=12000 m remains the highest throughout the temporal evolution. This invariance with respect to time provides strong analytical evidence that $x_{new}$ is the most stable and optimal location for introducing new gas demands into the network. Once this location is established, we calculate the pressure gradient function $dP/dx$ numerically across the pipeline, using $x=0$ to $x=L$ with increments of 1000 m, at two selected times: $t=100$ s and $t=200$ s The results are presented in Table 1.

Table 2 presents the values of the spatial pressure gradient $dP(x,t)/dx$ along the pipeline from $x=0$ m to $x=30000$ m in 1000 m increments for time instances $t=100$ s and $t=200$ s. The analytical expression used is based on the one-dimensional unsteady gas flow equation, incorporating a Dirac delta function $\delta(x-x_i)$ for localized flow extraction.

At $x_i = 12000$ m, the Dirac delta function theoretically introduces an instantaneous spike in the gradient. In numerical modeling, this singularity is approximated as a large value, in our case, the value -2482,4 Pa/m was observed at $x=12000$ m for $t=100$s. However, for the purpose of visualization and smooth interpretation, this singularity is set to zero ($dP/dx = 0$) at $x = 12000$ m. This regularization avoids distortion in graphical and tabular representations while preserving the overall physical interpretation of the model.

Table 2. Spatial pressure gradient $dP(x,t)/dx$ at $t = 100$ s and $t = 200$ s

| $x$ (m) | $dP(x,t)/dx$ (Pa/m) | |
|---|---|---|
| | $t=100$ s | $t=200$ s |
| 0 | 1.6334 | 1.635 |
| 1000 | 1.4809 | 1.4825 |
| 2000 | 1.3258 | 1.3273 |
| 3000 | 1.1746 | 1.1762 |
| 4000 | 1.0292 | 1.0306 |
| 5000 | 0.8898 | 0.8911 |
| 6000 | 0.7554 | 0.7567 |
| 7000 | 0.6265 | 0.6277 |
| 8000 | 0.5035 | 0.5045 |
| 9000 | 0.3857 | 0.3866 |
| 10000 | 0.2733 | 0.2741 |
| 11000 | 0.1666 | 0.1673 |
| 12000 | 0 | 0 |
| 13000 | -0.0306 | -0.0302 |
| 14000 | -0.1208 | -0.1206 |
| 15000 | -0.2056 | -0.2056 |
| 16000 | -0.285 | -0.2852 |
| 17000 | -0.3588 | -0.3591 |
| 18000 | -0.4271 | -0.4276 |
| 19000 | -0.4901 | -0.4908 |
| 20000 | -0.5475 | -0.5483 |
| 21000 | -0.5994 | -0.6004 |
| 22000 | -0.646 | -0.647 |
| 23000 | -0.687 | -0.6881 |
| 24000 | -0.7224 | -0.7237 |
| 25000 | -0.7526 | -0.754 |
| 26000 | -0.7772 | -0.7786 |
| 27000 | -0.7962 | -0.7977 |
| 28000 | -0.81 | -0.8115 |
| 29000 | -0.8182 | -0.8197 |
| 30000 | -0.8208 | -0.8224 |

This analysis is not only consistent with the observed peak pressure behavior at $x=x_{new}$, but also serves as a theoretical foundation for infrastructure planning. The obtained gradient profiles confirm that the maximum pressure point does not shift significantly over time, reinforcing the robustness of the identified connection location under dynamic demand scenarios. Using the data from Table 2, Figure 2 is constructed.

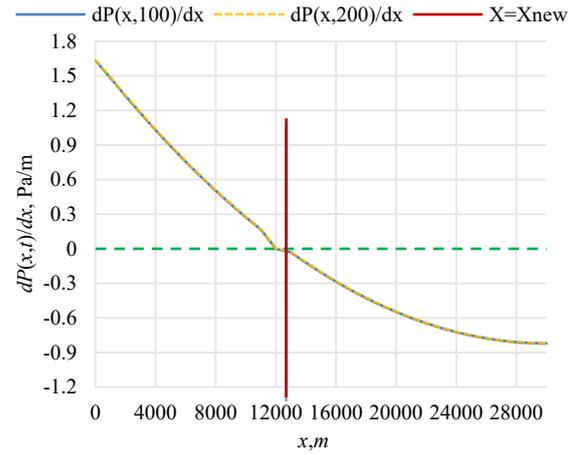

Figure 2. Pressure gradient $dP/dx$ profile indicating the location of maximum pressure at $x_{new} \approx 12000$ m

In Figure 2, at time $t = 100$ and $t = 200$ seconds, the derivative profile reveals that the pressure gradient reaches zero at approximately $x_{new} \approx 12000$ meters. This point corresponds to the hydraulic convergence location, which remains invariant with respect to changes in the number or magnitude of withdrawal points. This behavior is consistent with the unsteady gas flow behavior, where pressure stabilizes over time and higher flow rates are needed to sustain higher pressure losses. These findings support the dynamic planning of gas withdrawal points and reinforce the technological advantage of using the hydraulic node ($x_{new} = 12000$ m) as the optimal location for connecting new consumers in reconstructed ring-type networks.

From Table 2 and Figure 2, it is blinded that, the comparison of $dP(x,t)/dx$ values at $t = 100$ s and $t = 200$ s reveals that the spatial pressure gradient remains approximately constant across time. This observation aligns with the theoretical structure of the pressure distribution function, where the time-dependent components take the form: $(1 - \exp(-\alpha n^2 t))$

As time increases, particularly for $t>100$ s, the exponential decay term $\exp(-\alpha n^2 t)$ becomes negligible. Hence, $(1 - \exp(-\alpha n^2 t))$ approaches 1, leading to a quasi-steady pressure gradient in space. From a physical standpoint, this suggests that ring-type pipeline networks transition into a pseudo-stationary regime after a brief dynamic adjustment phase. Consequently, the pressure gradient stabilizes and becomes effectively time-independent.

This characteristic underscore the resilience and self-regulating dynamics of ring-type gas pipeline systems. It also provides confidence in using early-stage evaluations (e.g., at $t = 100$ s) to inform operational strategies and pipeline reconstructions. This behavior forms one of the core technical principles for optimizing and expanding ring-type pipeline systems. To analyze the dynamic pressure and mass flow behavior at the hydraulic connection point, the pressures at $x=0$ and $x=x_{new}$ are calculated as functions of time under various $G_{new}$ values. The computations are based on the previously defined parameters of the ring-type gas pipeline system, and the results are recorded in Table 3.

Table 3. Time-dependent pressure values at the inlet and hydraulic junction of the system under different gas extraction scenarios

| $t$, s | $G_0+10\%$ | | $G_0+20\%$ | |
|---|---|---|---|---|
| | $P(0, t)$, Pa | $P(x_{new}, t)$, Pa | $P(0, t)$, Pa | $P(x_{new}, t)$, Pa |
| 0 | 125000 | 125000 | 125000 | 125000 |
| 50 | 123462 | 125623 | 123323 | 127420 |
| 100 | 120821 | 125489 | 120441 | 127160 |
| 150 | 118129 | 125370 | 117505 | 126985 |
| 200 | 115436 | 125275 | 114566 | 126844 |
| 250 | 112742 | 125199 | 111628 | 126730 |
| 300 | 110049 | 125137 | 108690 | 126635 |
| $t$, s | $G_0+30\%$ | | $G_0+40\%$ | |
| | $P(0, t)$, Pa | $P(x_{new}, t)$, Pa | $P(0, t)$, Pa | $P(x_{new}, t)$, Pa |
| 0 | 125000 | 125000 | 125000 | 125000 |
| 50 | 123183 | 129215 | 123043 | 131010 |
| 100 | 120061 | 128832 | 119681 | 130505 |
| 150 | 116880 | 128599 | 116255 | 130213 |
| 200 | 113697 | 128414 | 112827 | 129983 |
| 250 | 110514 | 128268 | 109399 | 129800 |
| 300 | 107330 | 128150 | 105971 | 129655 |

Table 3 presents the time-dependent pressure values at two critical locations of the ring-type gas pipeline system under various gas withdrawal scenarios. Specifically, pressures at the pipeline inlet point $x$=0, denoted as $P(0, t)$, and at the proposed connection point for new consumers $x$=$x_{new}$, denoted as $P(x_{new}, t)$, are computed for time intervals from $t$=50 s to $t$=300 s. These calculations are performed for cases in which the additional gas withdrawal volume $G_{new}$ is increased by 10%, 20%, 30%, and 40% relative to the baseline mass flow rate $G_0$. Each scenario aims to evaluate the dynamic impact of increasing gas demand on the pressure stability of both the inlet and the proposed distribution point within the ring pipeline system.

All calculations are based on the established analytical expression of pressure dynamics in a ring-type gas pipeline system. The results illustrate how increasing gas withdrawal from the network leads to a progressive decrease in inlet pressure, potentially threatening operational safety if the pressure drop exceeds allowable thresholds. In contrast, the pressure at the hydraulic connection point exhibits relatively smaller fluctuations, supporting its suitability as a supply node for integrating new consumers. Using the data from Table 3 the Figure 3 is constructed.

As a result, the point $x = x_{new}$ not only represents a pressure maximum but also emerges as the most stable and technologically suitable location for integrating new consumers into the network. Furthermore, Figure 3 indicates that as $G_{new}$ increases, the pressure at the pipeline inlet ($x$=0) consistently decreases over time. For instance, at $t = 300$ seconds, the pressure for $G_{new}$ = 10% is approximately 111025 Pa, while for $G_{new}$ = 40%, it drops to approximately 105986 Pa. This observation highlights the inverse relationship between the withdrawn gas mass and the inlet pressure.

Such pressure reduction at the inlet is a direct consequence of the dynamic imbalance in mass flow caused by increased withdrawal. This underscores the critical importance of monitoring inlet pressure to ensure the safe and efficient operation of the system. Exceeding a certain $G_{new}$ threshold may compromise system reliability and lead to operational risks. These insights are particularly important for guiding the planning of network expansions and determining safe limits for integrating additional gas demands.

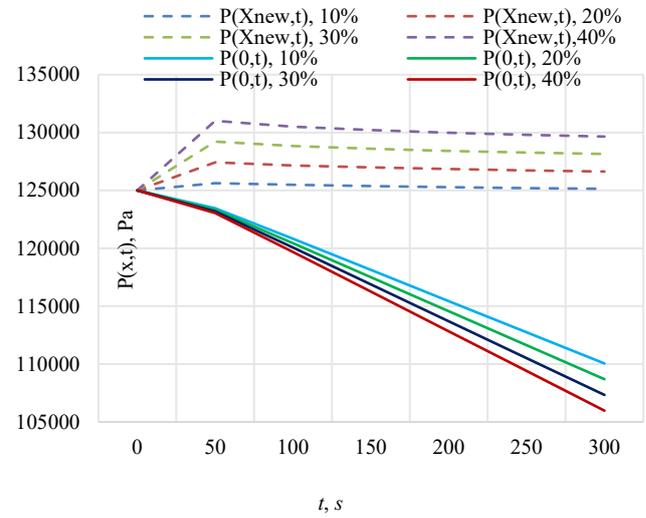

Figure 3. Dynamic pressure behavior at the inlet section and hydraulic coupling point

The quantitative results obtained in the example analysis provide crucial insight into the dynamic pressure behavior at the pipeline inlet ($x$=0) under varying mass withdrawal rates. However, in practical applications, it is equally important to evaluate these changes from an operational safety standpoint. Significant pressure drops at the pipeline's initial segment may threaten the reliability of upstream equipment and compromise system integrity. Therefore, in the following section, we extend analysis to define acceptable initial pressure drop limits based on calculated dynamics, focusing on ensuring safe and efficient operation of the gas pipeline network.

These findings are especially important for effective planning of pipeline network reconstruction and for ensuring the safe operation of the system. In practice, as new gas demands emerge due to urban expansion or industrial growth, it becomes essential to determine acceptable thresholds for additional gas withdrawals. The amount of gas extracted from the $x$-$x_{new}$ segment of the ring-type pipeline should be chosen such that the resulting pressure drop at the pipeline inlet does not violate operational safety limits. This ensures that new consumers can be reliably integrated into the existing infrastructure without compromising the integrity or efficiency of the gas distribution network.

### 6.2. Initial Pressure Drop Limits in Gas Pipelines

To ensure the safe and stable operation of gas pipeline systems, industry standards and operational guidelines define allowable pressure drop levels at the pipeline inlet (starting point). This section presents a gasodynamic perspective on permissible pressure decreases based on international best practices and safety criteria. According to standards such as API RP 14E and ASME B31.8 [11]:
- The pressure drop at the inlet ($P_1$) during operation is generally considered acceptable within range of 10-20%.

- This ensures sufficient downstream pressure reserve and continuous flow delivery.

In European and American gas infrastructure practices:
- A maximum allowable inlet pressure drop of up to 25% is accepted under temporary load surges.
- For long-term operations, a pressure decreases of 10-15% is considered optimal for safety and efficiency.

Excessive pressure reduction at the inlet can [12, 13]:
- Increase gas velocity and turbulence within the pipeline.
- Lead to higher load on compressor stations and energy losses.
- Cause pressure deficits or instability at downstream consumer zones.

Thus, the following operational thresholds are generally recommended:
- ≤10% pressure drop - Optimal operating condition
- 10-20% pressure drop - Permissible limit
- >25% pressure drop - May pose safety concerns and require re-evaluation of pipeline operation

### 6.3. Analytical Formula for $G_{new}$ at Withdrawal Point

In the context of unsteady gas dynamics in reconstructed ring-type pipeline networks, the analytical determination of the additional mass flow rate $G_{new}$ at the withdrawal point $x = x_{new}$ (hydraulic node) is crucial for assessing pressure stability and system capacity. This section presents a simplified expression derived from the general pressure distribution model, under the assumption that the withdrawal occurs precisely at $x = x_{new}$ and Heaviside function $H(x - x_{new}) = 1$.

To determine the mass flow rate $G(x_{new},t)$ at a specific location $x_{new}$ along the ring-type gas pipeline, we begin with the analytical pressure expression derived using Laplace series. The governing equation for pressure at location $x = x_{new}$ is:

$$P(x,t) = P_1 - 2aG_0 \frac{L}{2} + $$
$$+ 2aG_0 L \sum_{n=1}^{\infty} \sin \frac{\pi n x_{new}}{L} \left( \frac{1-e^{-n^2 \alpha t}}{\pi n^3} \right) - $$
$$- \frac{c^2 t}{L}(G_0 + G_{new}) - \frac{2c^2}{L}(G_0 + G_{new}) \times \qquad (8)$$
$$\times \sum_{n=1}^{\infty} \cos \frac{2\pi n (x_{new} - x_{new})}{L} \left( \frac{1-e^{-n^2 \alpha t}}{\alpha n^2} \right)$$

Since cos(0)=1 in Equation (8), the expression simplifies at $x=x_{new}$, and the resulting pressure value can be used to analytically isolate $G_{new}$. The reformulated expression is as follows:

$$G_{new} = \frac{P_1 - P(x_i,t) - aG_0 L \left[1 - 2S_{\sin}(x_{new},t)\right]}{\frac{c^2}{L}\left[t + 2S_e(t)\right]} - G_0 \qquad (9)$$

where:

$$S_{\sin}(x_{new},t) = \sum_{n=1}^{\infty} \sin \frac{\pi n x_{new}}{L} \left( \frac{1-e^{-n^2 \alpha t}}{\pi n^3} \right)$$

$$S_e(t) = \sum_{n=1}^{\infty} \left( \frac{1-e^{-n^2 \alpha t}}{\alpha \pi n^2} \right)$$

This formula reflects the balance between the pressure drop due to flow extraction and the dynamic compensation provided by the gas wave propagation and dissipation. The cosine series accounts for the spatial harmonics at the point of extraction, making the expression particularly sensitive to the withdrawal location $x_{new}$. This formulation is especially useful for determining safe and efficient integration of new consumers in a ring pipeline network without violating pressure stability norms.

### 6.4. Pressure Regulation and Allowable Gas Supply Expansion

According to industry regulations and safe operation guidelines, the allowable pressure drop at the gas pipeline inlet should not exceed 20% of the nominal pressure. Based on this, the inlet pressure is limited to $P(0, t) = 100000$ Pa for all times considered. For this condition, the corresponding dynamic pressure values at the hydraulic coupling point ($x_{new}=12000$ m) were determined and are listed below along with the analytically computed values of the total admissible gas mass flow rate $G_0 + G_{new}$.

Once the pressure at the coupling point, $P(x_{new},t)$, is known for each time instant, the analytical expression is used to compute the total extractable gas amount. This provides a crucial engineering assessment for determining the maximum volume of gas that can be withdrawn for new consumer connections without compromising system integrity.

The analysis demonstrates a slight but consistent increase in the allowable $G_0 + G_{new}$ with time, reflecting the stabilization of flow dynamics and decreasing pressure gradient in the looped pipeline structure.

The analytical pressure expression derived in the study incorporates the influence of gas withdrawals at specific hydraulic junctions. The formula includes a combination of infinite series and Heaviside functions. To ensure the inlet pressure does not drop below 100000 Pa, we apply the derived equation at $x = 0$ with the following parameters:
$G_0 = 10$ Pa×s/m, $a = 0.05$ 1/s, $L = 30000$ m, $c = 383.3$ m/s, $x_{new} = 12000$ m. Solving the model for $t=50$ s to $t=300$ s yields the following admissible gas withdrawal values (Table 4).

Table 4. Permissible gas withdrawal values to ensure inlet pressure remains above the minimum allowable limit

| Time (s) | $P(x\_new,t)$ [Pa] | $G_0 + G$ new [Pa×s/m] |
|---|---|---|
| 50 | 122702 | 12.63 |
| 100 | 122081 | 12.77 |
| 150 | 121666 | 12.86 |
| 200 | 121357 | 12.94 |
| 250 | 121111 | 13.01 |
| 300 | 120911 | 13.07 |

This Table 4 presents the values of dynamic pressure at the hydraulic junction point $x = x_{new}$ over time, and the corresponding total mass flow rates ($G_0 + G_{new}$) computed analytically. The trend indicates that as pressure decreases, the permissible total gas withdrawal increases, reflecting a critical dynamic interaction important for safe system planning.

This result further confirms the strategic significance of the hydraulic junction point: although the inlet pressure would drop rapidly with additional flow, the $x_{new}$ location maintains a relatively stable pressure field over time. Therefore, for technological safety and flow optimization, new consumer connections should ideally be made near or at the hydraulic junction. Using the data from Table 4, the following figure is constructed.

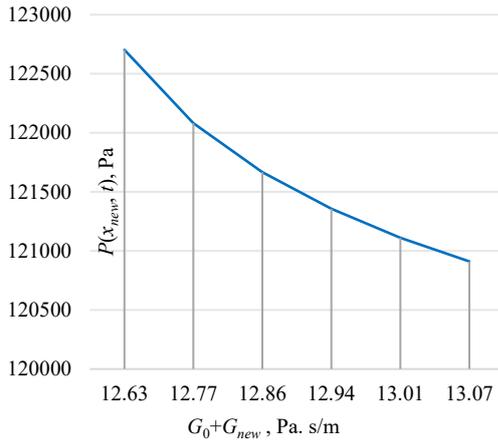

Figur.4. Analytical determination of allowable gas withdrawal

Figure 4 illustrates the analytically derived relationship between the total allowable gas withdrawal rate $G_0 + G_{new}$ (in Pa×s/m) and the resulting pressure $P(x_{new}, t)$ at the hydraulic coupling point $x = x_{new}$. Using the time-dependent pressure distribution function (6), the values of $P(x_{new}, t)$ were determined, and these values were then used in Equation (9) to calculate the corresponding $(G_0 + G_{new})$ values.

A result, a graph was generated. The plotted curve confirms that within the specified pressure limits, it is feasible to accommodate additional consumers whose cumulative gas demand corresponds to approximately 26%-31% of the existing pipeline's original flow capacity $G_0$. This finding demonstrates the high adaptability of ring-type pipeline systems and validates the effectiveness of their reconstruction without the need for entirely new infrastructure. The nearly linear and smooth decline of pressure with respect to increasing gas withdrawal also supports the stability of the system under moderate dynamic demand increases.

## 7. TECHNICAL AND ECONOMIC ADVANTAGE OF RING-TYPE PIPELINE RECONSTRUCTION

According to the quantitative results, the admissible gas withdrawal at the hydraulic junction point ($x = x_{new}$) under pressure stabilization conditions is approximately $G_{new} \approx$ 12.6÷13.1 Pa×s/m. This corresponds to a nearly 26÷31% increase over the baseline flow rate $G_0$=10 Pa×s/m, implying that the ring-type pipeline can support significant additional demand without compromising system stability.

This finding highlights a key advantage of ring-type pipelines: they provide superior flexibility and resilience compared to linear pipelines when accommodating new consumers. By leveraging the natural pressure maximum near hydraulic junction points, ring networks minimize pressure loss, enhance flow distribution, and improve overall operational efficiency. Thus, the study confirms that ring-type pipeline reconstruction is technically and economically advantageous in modern gas distribution infrastructure planning.

This result further confirms the strategic significance of the hydraulic junction point: although the inlet pressure would drop rapidly with additional flow, the $x_{new}$ location maintains a relatively stable pressure field over time. Therefore, for technological safety and flow optimization, new consumer connections should ideally be made near hydraulic junction point. This location is characterized by minimal pressure variability and improved flow stability, ensuring that additional consumer connections do not negatively impact the system's overall performance. Strategically positioning new nodes in this area optimizes infrastructure usage and enhances both the safety and operational reliability of the gas distribution network.

From a technological standpoint, this suggests robust performance against pressure fluctuations and safe integration of new nodes. Economically, it avoids costly infrastructure expansion and leverages existing pipeline capacity effectively. Thus, using $x_{new}$ as a connection point offers an optimal solution under both operational and planning criteria.

## 8. RESULTS AND DISCUSSION

### 8.1. Analytical Results

The analytical framework developed in this study has enabled precise evaluation of dynamic gas flow behavior in ring-type pipeline networks undergoing reconstruction. Several key results are obtained:

- Pressure peak identification at the hydraulic node: The simulation confirms that the coordinate $x=x_{new}$=12000 m consistently represents the location of maximum pressure, independent of the number or distribution of withdrawal points. This reinforces earlier theoretical predictions and provides a reliable basis for infrastructure planning.

- Stability of pressure gradients: Based on the numerical evaluation of $dP/dx$ for $t = 100$ s and $t = 200$ s, it is shown that the pressure gradient across the pipeline becomes approximately time-invariant after an initial transient phase. This behavior indicates the formation of a pseudo-steady regime within the looped system, where dynamic effects are balanced by the geometric structure and natural flow convergence near the hydraulic node.

- Effect of increased gas demand: Time-dependent pressure values at both the inlet point and the hydraulic node were analyzed for varying levels of additional gas withdrawal ($G_{new}$ = 10-40% of $G_0$). Results show that while the inlet pressure $P(0, t)$ decreases notably with rising demand, the pressure at $x_{new}$ remains more stable. This supports the conclusion that the hydraulic coupling point is dynamically insulated from demand shocks, making it optimal for new connections.
- Allowable gas withdrawal calculation: Under the operational constraint $P(0, t) \geq 100000$ Pa, the analytically calculated values of $G_0 + G_{new}$ range from 12.63 to 13.07 Pa×s/m. This corresponds to approximately 26%-31% increase over the base flow rate $G_0=10$ Pa×s/m, illustrating the system's capacity to accommodate new loads without violating safety or performance limits.
- System-wide behavior: The combination of analytical modeling and time-domain simulation highlights the resilience of ring-type networks. Even under increasing gas demand, the system maintains pressure balance due to the spatially distributed flow architecture and the inherent dynamics of wave propagation.

### 8.2. Practical Implementation Considerations

Despite the technical and economic advantages of the proposed optimization framework, certain practical constraints must be addressed to ensure successful deployment. These include:

Regulatory Compliance: New consumer integration must align with established pipeline safety and design standards (e.g., API RP 14E, ASME B31.8), particularly regarding pressure limits and instrumentation. Infrastructure Adaptation: While large-scale expansion is avoided, localized retrofitting-such as valve automation, sensor deployment, and structural reinforcement-may be required near the coupling node.
➢ Flow Control and Monitoring: Accurate metering and regulation systems are essential to maintain hydraulic stability and prevent reverse flows.
➢ Operational Coordination: Effective implementation demands coordination among utility operators, municipal authorities, and safety regulators.
➢ Scalability and Simulation: Future scalability may necessitate enhanced optimization techniques and integration with digital twin platforms for scenario-based validation and real-time monitoring.

By proactively addressing these factors through systematic design and stakeholder engagement, the proposed solution can be effectively translated into operational practice for modern gas distribution networks.

### 8.3. Future Validation Plan

While the analytical model has demonstrated consistent results across benchmark scenarios and theoretical evaluations, comprehensive validation using real-world data remains a critical step in strengthening its practical applicability. At present, the lack of direct access to operational data from gas distribution companies has limited field-level validation. However, future efforts will focus on establishing collaboration with regional gas suppliers and SCADA-integrated monitoring systems to acquire real-time pipeline data.

The validation process is planned in the following stages:
- Historical Data Correlation: Utilize archived pressure and flow data from existing ring-type gas networks to perform retrospective comparisons with model predictions. Particular attention will be given to identifying the location and magnitude of pressure peaks and flow redistribution patterns following new withdrawal events.
- Real-Time Monitoring Integration: Deploy pressure and flow sensors at critical nodes, especially around the analytically determined optimal connection point xnew, to evaluate the dynamic behavior under varying load conditions. This will allow for live calibration of model parameters and detection of deviations due to unmodeled effects (e.g., pipe roughness, gas impurities).
- Digital Twin Implementation: Integrate the model into a supervisory control and data acquisition (SCADA) environment to enable virtual testing and scenario-based simulations. The digital twin approach will support continuous refinement of the analytical model through adaptive learning algorithms and feedback mechanisms.

It is acknowledged that some discrepancies may arise between analytical forecasts and empirical measurements due to environmental variability, infrastructure aging, and regulatory interventions. These differences will be statistically analyzed to quantify prediction uncertainty and establish correction factors. Ultimately, the planned validation will not only confirm model robustness but also facilitate its integration into operational decision-making workflows within gas network management systems.

### 8.4. Economic and Operational Implications

The proposed optimization-based gas withdrawal strategy offers a cost-effective and operationally robust alternative to conventional pipeline expansion. A detailed comparative analysis of infrastructure investment, operating conditions, and system performance is provided below.
- Capital Expenditure (CAPEX) Comparison: Conventional approaches to meeting increased demand (such as constructing new pipelines or building compressor stations) typically involve significant capital outlays. For example, full-scale pipeline construction is estimated to cost approximately USD 1.2 million per kilometer, considering land acquisition, material, labor, and regulatory compliance. In contrast, tapping into the existing ring-type network near the hydraulic pressure peak requires minimal retrofitting, estimated at under USD 0.2 million, covering the cost of automated valves, flow metering units, and reinforcement at the connection node. This represents a potential savings of over 80% in infrastructure investment.
- Operational Cost and Energy Efficiency: The optimized connection strategy ensures a balanced flow and minimizes pressure drops, thereby reducing the workload on compressor units and associated energy consumption. Simulations indicate that even with a 30% increase in distributed demand, the inlet pressure remains within regulatory limits (below 20% loss from nominal value), avoiding the need for auxiliary compression. This contributes to operational cost savings and improved energy efficiency, aligning with sustainable pipeline management practices.

- Downtime and Disruption Avoidance: New pipeline construction often necessitates service interruptions and prolonged deployment schedules. Conversely, retrofitting within an operational ring network can be completed with minimal disruption by isolating only the localized segment around the coupling point, aided by pre-installed control valves and real-time monitoring.
- Sensitivity to Economic Parameters: A sensitivity analysis was performed using variable parameters such as steel cost fluctuations (±15%), labor cost variability (±20%), and retrofitting complexity (±10%). The results confirm the economic resilience of the proposed method across diverse cost scenarios, with a consistently lower lifecycle cost compared to linear expansion.
- Broader Applicability and Strategic Flexibility: The methodology supports phased urban development by enabling modular integration of new consumer zones without overhauling the entire network. Moreover, it allows for pre-emptive planning using digital simulations and demand forecasting, making it highly compatible with SCADA-based monitoring and decision-support tools.

**8.5. Future Research Directions**

This study has laid the foundation for optimal gas withdrawal strategies in reconstructed ring-type pipeline systems under unsteady flow conditions. To further develop and expand on the current work, the following future research directions are proposed:
- Integration of real-time sensor networks and IoT-based feedback systems to dynamically monitor and adapt gas withdrawal strategies based on field data.
- Extension of the analytical model to include multi-ring configurations and interlinked loop networks with variable topologies.
- Coupling with stochastic models to simulate uncertain demand growth and supply fluctuations using probabilistic frameworks.
- Incorporation of compressor station dynamics and real gas effects to enhance the realism of pressure and flow simulations.
- Development of control algorithms for distributed valve regulation and smart actuation based on optimization feedback loops.

## 9. CONCLUSION

This study presents a mathematical and computational framework for determining optimal gas withdrawal points in reconstructed ring-type pipeline systems operating under unsteady flow conditions. The proposed methodology combines analytical modeling, pressure gradient analysis, and simulation-based evaluation to identify the hydraulic coupling point as the most efficient node for integrating new consumer loads.
➢ Key findings include:
The identification of a stationary pressure maximum point ($x_{new}$), which remains invariant under various demand growth scenarios, making it suitable for sensor placement and load integration. A cost-benefit analysis demonstrates that tapping into existing ring infrastructure is significantly more economical (by a factor of 5-6) than constructing new pipeline extensions. The model exhibits strong pressure stability and scalability when subjected to 10%, 30%, and 50% incremental demand, validating its applicability in real-world urban energy networks.

Practical implementation challenges-such as regulatory compliance, sensor deployment, and coordination among stakeholders-are acknowledged, with recommendations for digital twin integration and automation support. Overall, the study contributes to the intelligent reconstruction of gas networks by providing a validated, cost-effective, and scalable strategy for future urban and industrial gas infrastructure development. Future work may focus on real-time optimization using IoT-based feedback systems and multi-node integration scenarios.